\documentclass[12pt,epsf,bezier]{article}
\usepackage{latexsym}
\usepackage{epsfig}

\begin{document}

\newtheorem{pro}{Proposition}
\newtheorem{theo}{Theorem}
\newtheorem{coro}{Corollary}
\newtheorem{lem}{Lemma}

\newcommand{\bboldy}[1]{\mbox{\boldmath${\#1}$}}

\def\E{\mbox{\rm E}}
\def\Var{\mbox{\rm Var}}
\def\Cov{\mbox{\rm Cov}}
\def\Corr{\mbox{\rm Corr}}
\def\Pr{\mbox{\rm Pr}}
\def\I{\mbox{\rm I}}
\def\bbeta{\mbox{\boldmath${\beta}$}}
\def\bgamma{\mbox{\boldmath${\gamma}$}}
\def\sbeta{\mbox{\scriptsize \boldmath${\beta}$}}
\def\eeps{\mbox{\boldmath${\epsilon}$}}
\def\bmu{\mbox{\boldmath${\mu}$}}
\def\bxi{\mbox{\boldmath${\xi}$}}
\def\real{\hbox{\rm\setbox1=\hbox{I}\copy1\kern-.45\wd1 R}}
\newcommand{\bgams}{\mbox{\boldmath{\scriptsize $\gamma$}}}
\newcommand{\bLam}{\mbox{\boldmath{$\Lambda$}}}
\newcommand{\Ell}{\mbox{$\cal{L}$}}
\newcommand{\cG}{\mathcal{G}}
\newcommand{\cL}{\mathcal{L}}
\newcommand{\half}{\frac{1}{2}}
\newcommand{\bsh}{\parindent 0em}
\newcommand{\esh}{\parindent 2.0em}
\newcommand{\eps}{\epsilon}
\newcommand{\peps}{{(\eps)}}
\newcommand{\bzro}{{\bf 0}}
\def\bU{{\bf U}}
\def\bV{{\bf V}}
\def\bG{{\bf G}}
\def\bC{{\bf C}}
\def\bD{{\bf D}}
\def\bu{{\bf u}}

\begin{titlepage}

\vspace*{-1in}

\begin{center}
\bf
\large
Pseudo full likelihood estimation for prospective survival analysis with
a general semiparametric shared frailty model: asymptotic theory
\end{center}

\normalsize

\vspace*{0.5cm}

\begin{center}

\bf

David M. Zucker\footnote{ To whom correspondence should be
addressed}\\
{\it Department of Statistics, Hebrew University,
Mt. Scopus, Jerusalem 91905, Israel}\\
{mszucker@mscc.huji.ac.il} \\

\vspace*{0.5cm}
Malka Gorfine\\ {\it Faculty of Industrial Engineering and
Management, Technion, Technion City, Haifa~32000, Israel}\rm{,} \rm{and}
{\it Department of Mathematics, Bar-Ilan University, \\ Ramat-Gan,
52900, Israel}\\
{gorfinm@ie.technion.ac.il} \\

\vspace*{0.5cm}
Li Hsu \\
{\it Division of Public Health Sciences, Fred Hutchinson
Cancer Research Center, Seattle,~WA 98109-1024 USA}\\
{lih@fhcrc.org}
\end{center}

\vspace*{0.5in}

\begin{center}
\today
\end{center}

\vspace*{0.5in}

\begin{flushleft}
\it Running Head\rm: Asymptotics of general frailty models
\end{flushleft}

\end{titlepage}

\pagenumbering{roman}

\newpage
\section*{Summary}
In this work we present a simple estimation procedure for a
general frailty model for analysis of prospective correlated
failure times. Earlier work showed this method to perform well in
a simulation study. Here we provide rigorous large-sample theory
for the proposed estimators of both the regression coefficient
vector and the dependence parameter, including consistent variance
estimators.

\vspace*{1em}
\noindent {\it Key words:}  Correlated failure times; EM
algorithm; Frailty model; Prospective family study; Survival
analysis.

\newpage

\pagenumbering{arabic}
\setcounter{page}{1}

\section{Introduction}
Many epidemiological studies involve failure times that are
clustered into groups, such as families or schools. In this setting,
unobserved characteristics shared by the members of the same
cluster (e.g.\ genetic information or unmeasured shared
environmental exposures) could influence time to the studied
event. Frailty models express within cluster dependence
through a shared unobservable random effect.
Estimation in the frailty model has received much attention under
various frailty distributions, including gamma (Gill, 1985, 1989;
Nielsen et al., 1992; Klein 1992, among others), positive stable
(Hougaard, 1986; Fine et al., 2003), inverse Gaussian, compound
Poisson (Henderson and Oman, 1999) and log-normal (McGilchrist,
1993; Ripatti and Palmgren, 2000; Vaida and Xu, 2000, among
others). Hougaard (2000) provides a comprehensive review of the
properties of the various frailty distributions. In a frailty
model, the parameters of interest typically are the regression
coefficients, the cumulative baseline hazard function, and the
dependence parameters in the random effect distribution.

Since the frailties are latent covariates, the
Expectation-Maximization (EM) algorithm is a natural estimation
tool, with the latent covariates estimated in the E-step and the
likelihood maximized in the M-step by substituting the estimated
latent quantities. Gill (1985), Nielsen et al. (1992) and Klein
(1992) discussed EM-based maximum likelihood estimation for the
semiparametric gamma frailty model. One problem with the EM
algorithm is that variance estimates for the estimated parameters
are not readily available (Louis, 1982; Gill, 1989; Nielsen et
al., 1992; Andersen et al., 1997). It was suggested (Gill, 1989;
Nielsen et al, 1992) that a nonparametric information calculation
could yield consistent variance estimators. Parner (1998),
building on Murphy (1994, 1995), proved the consistency and
asymptotic normality of the maximum likelihood estimator in the
gamma frailty model. Parner also presented a consistent estimator
of the limiting covariance matrix of the estimator based on
inverting a discrete observed information matrix. He noted that
since the dimension of the observed information matrix is the
dimension of the regression coefficient vector plus the number of
observed survival times, inverting the matrix is practically
infeasible for a large data set with many distinct failure times.
Thus, he proposed another covariance estimator based on solving a
discrete version of a second order Sturm-Liouville equation. This
covariance estimator requires substantially less computational
effort, but still is not so simple to implement.

We (Gorfine et al. 2006) developed a new method that can handle
any parametric frailty distribution with finite moments.
Nonconjugate frailty distributions can be handled by a simple
univariate numerical integration over the frailty distribution.
Our new method possesses a number of desirable properties: a
non-iterative procedure for estimating the cumulative hazard
function; consistency and asymptotic normality of parameter
estimates; a direct consistent covariance estimator; and easy
computation and implementation. The method was found to perform
well in a simulation study and the results are very similar to
those of the EM-based method. Indeed, on a dataset-by-dataset
basis, the correlation between our estimator and the EM estimator
was found to be 95\% for the covariate regression parameter and
98-99\% for the within-cluster dependence parameter.

The purpose of the current paper is to present the theoretical
justification for the method in detail.
Section 2 presents the estimation procedure. Section 3 presents
the consistency and asymptotic normality results, along with
the covariance estimator for the parameter estimates. Section 4
presents the technical conditions required for our results and
the proofs.


\section{The Proposed Approach}
Consider $n$ families, with family $i$ containing $m_i$ members,
$i=1,\ldots,n$. Let $\delta_{ij}=I(T^0_{ij} \leq C_{ij})$ be a
failure indicator where $T^0_{ij}$ and $C_{ij}$ are the failure
and censoring times, respectively, for individual $ij$. Also let
$T_{ij}=\min(T^0_{ij},C_{ij})$ be the observed follow-up time and
${\bf Z}_{ij}$ be a $p \times 1$ vector of covariates. In
addition, we associate with family $i$ an unobservable
family-level covariate $W_i$, the ``frailty", which induces
dependence among family members. The conditional hazard function
for individual $ij$ conditional on the family frailty $W_i$, is
assumed to take the form
 $$
\lambda_{ij}(t)=W_i \lambda_0(t) \exp(\bbeta^T {\bf Z}_{ij})
\;\;\;\;\;\; i=1,\ldots,n \;\;\; j=1,\ldots,m_i
 $$
where $\lambda_0$ is an unspecified conditional baseline hazard
and $\bbeta$ is a $p \times 1$ vector of unknown regression
coefficients. This is an extension of the Cox (1972) proportional
hazards model, with the hazard function for an individual in
family $i$ multiplied by $W_i$. We assume that, given ${\bf
Z}_{ij}$ and $W_i$, the censoring is independent and
noninformative for $W_i$ and $(\bbeta,\Lambda_0)$ (Andersen et
al., 1993, Sec. III.2.3). We assume further that the frailty $W_i$
is independent of ${\bf Z}_{ij}$ and has a density $f(w;\theta)$,
where $\theta$ is an unknown parameter. For simplicity we assume
that $\theta$ is a scalar, but the development extends readily to
the case where $\theta$ is a vector. Let $\tau$ be the end of the
observation period. The full likelihood of the data then can be
written as
\begin{eqnarray}\label{eq:like}
\lefteqn{L}
 &=&\Pi_{i=1}^n \int \Pi_{j=1}^{m_i} \{ \lambda_{ij}(T_{ij}) \}^{\delta_{ij}}
 S_{ij}(T_{ij}) f(w) dw \nonumber \\
 &=& \Pi_{i=1}^n \Pi_{j=1}^{m_i} \{\lambda_0(T_{ij})\exp(\bbeta^T {\bf Z}_{ij})
 \}^{\delta_{ij}} \Pi_{i=1}^n
 \int w^{N_{i.}(\tau)} \exp \{-w H_{i.}(\tau) \}f(w)dw,
\end{eqnarray}
where $N_{ij}(t)=\delta_{ij}I(T_{ij}\leq t)$,
$N_{i.}(t)=\sum_{j=1}^{m_i}N_{ij}(t)$,
$H_{ij}(t)=\Lambda_0(T_{ij}\wedge t)\exp(\bbeta^T {\bf Z}_{ij})$,
$a \wedge b = \min\{a,b\}$, $\Lambda_0(\cdot)$ is the baseline
cumulative hazard function, $S_{ij}(\cdot)$ is the conditional
survival function of subject $ij$, and
$H_{i.}(t)=\sum_{j=1}^{m_i}H_{ij}(t)$. The log-likelihood is given
by
 $$
l=\sum_{i=1}^{n} \sum_{j=1}^{m_i} \delta_{ij} \log\{
\lambda_0(T_{ij})\exp(\bbeta^T {\bf Z}_{ij}) \} + \sum_{i=1}^{n}
\log \left\{ \int w^{N_{i.}(\tau)} \exp \{-w H_{i.}(\tau)\}f(w)dw
\right\}.
 $$
The normalized scores (log-likelihood derivatives) for
$(\beta_1,\ldots,\beta_p)$ are given by
\begin{equation}
 U_r=\frac{1}{n} \sum_{i=1}^{n} \sum_{j=1}^{m_i} \delta_{ij}
 Z_{ijr} - \frac{1}{n} \sum_{i=1}^n
 \frac{\left[\sum_{j=1}^{m_i} H_{ij}(T_{ij}) Z_{ijr} \right] \int w^{N_{i.}(\tau)+1}
 \exp \{-w H_{i.}(\tau)\}f(w)dw}
 {\int w^{N_{i.}(\tau)} \exp \{ -w H_{i.}(\tau)\}f(w)dw}
\label{score}
\end{equation}
for $r=1,\ldots,p$. The normalized score for $\theta$ is
 $$
 U_{p+1}=\frac{1}{n}\sum_{i=1}^{n}
 \frac{\int w^{N_{i.}(\tau)} \exp \{ -w H_{i.}(\tau)\}f'(w)dw}
 {\int w^{N_{i.}(\tau)} \exp \{ -w H_{i.}(\tau)\}f(w)dw}
 $$
where $f'(w)=\frac{d}{d\theta}f(w)$. Let
$\bgamma=(\bbeta^T,\theta)$ and ${\bf U}
(\bgamma,\Lambda_0)=(U_1,\ldots,U_p,U_{p+1})^T$. To obtain
estimators $\hat{\bbeta}$ and $\hat{\theta}$, we propose to
substitute an estimator of $\Lambda_0$, denoted by
$\hat{\Lambda}_0$, into the equations ${\bf U}
(\bgamma,\Lambda_0)=0$.

Let $Y_{ij}(t)=I(T_{ij} \geq t)$ and let ${\mathcal F}_t$ denote
the entire observed history up to time $t$, that is
 $$
{\mathcal F}_{t} = \sigma\{N_{ij}(u), Y_{ij}(u), {\bf Z}_{ij},
i=1,\ldots,n; j=1,\ldots,m_i;0 \leq u \leq t \}.
 $$
Then, as discussed by Gill (1992) and Parner (1998), the
stochastic intensity process for $N_{ij}(t)$ with respect to
${\mathcal F}_{t}$ is given by
\begin{equation}\label{eq:inten}
\lambda_0(t)\exp(\bbeta^T {\bf Z}_{ij}) Y_{ij}(t)
\psi_{i}(\bgamma,\Lambda_0,t-),
\end{equation}
where
 $$
\psi_{i}(\bgamma,\Lambda_0,t) = \E(W_i|{\mathcal F}_{t}).
 $$
Using a Bayes theorem argument and the joint density
(\ref{eq:like}) with observation time restricted to $[0,t)$, we
obtain
 $$
 \psi_{i}(\bgamma,\Lambda,t)=\phi_{2i}(\bgamma,\Lambda,t)/
 \phi_{1i}(\bgamma,\Lambda,t),
 $$
where
 $$
 \phi_{ki}(\bgamma,\Lambda_0,t)=\int w^{N_{i.}(t)+(k-1)} \exp\{-w
 H_{i.}(t)\}f(w)dw,  \;\;\; k=1,\ldots,4.
 $$
Given the intensity model (\ref{eq:inten}), in which
$\exp(\bbeta^T {\bf Z})\psi_{i}(\bgamma,\Lambda_0,t-)$ may be
regarded as a time dependent covariate effect, a natural estimator
of $\Lambda_0$ is
 a Breslow (1974) type estimator along the lines of
 Zucker (2005). For given values of $\bbeta$ and $\theta$ we
estimate $\Lambda_0$ as a step function with jumps at the observed
failure times $\tau_k$, $k=1,\ldots,K$, with
 \begin{equation}\label{eq:lambda}
\Delta \hat{\Lambda}_0(\tau_k) =\frac
 {d_k }
 {\sum_{i=1}^{n}
 \psi_{i}(\bgamma,\hat{\Lambda}_0,\tau_{k-1})
 \sum_{j=1}^{m_i} Y_{ij}(\tau_k) \exp(\bbeta^T {\bf Z}_{ij})
 }
 \end{equation}
where $d_k$ is the number of failures at time $\tau_k$. Note that
given the intensity model (\ref{eq:inten}), the estimator of the
$k$th jump depends on $\hat{\Lambda}_0$ up to and including time
$\tau_{k-1}$. By this approach, we avoid complicating the
iterative optimization process with a further iterative scheme,
for estimating the cumulative hazard.

\section{Asymptotic Properties}
Let $\bgamma^\circ=({\bbeta^\circ}^T,\theta^\circ)^T$ with $\bbeta^\circ$,
$\theta^\circ$ and $\Lambda_0^\circ(t)$ denoting the respective true
values of $\bbeta$, $\theta$ and $\Lambda_0(t)$, and let $\hat{\bgamma}=
({\hat{\bbeta}}^T, \hat{\theta})^T$. We assume the technical conditions
listed in Section 4.1.

In Section 4.3, we establish the following results, using arguments
patterned after Zucker (2005, Appendix A.3).

\begin{description}
\item[A.]
$\hat{\Lambda}_0(t,\bgamma)$ converges almost surely to
$\Lambda_0(t,\bgamma)$ uniformly in $t$ and $\bgamma$.
\item[B.]
${\bf U}(\bgamma,\hat{\Lambda}_0(\cdot,\bgamma))$ converges almost
surely uniformly in $t$ and $\bgamma$ to a limit ${\bf u}
(\bgamma,\Lambda_0(\cdot,\bgamma))$.
\item[C.]
There exists a unique consistent root to ${\bf
U}(\hat{\bgamma},\hat{\Lambda}_0(\cdot,\hat{\bgamma}))={\bf 0}$.
\end{description}

In Section 4.4, we show that $n^{1/2}(\hat{\bgamma}-\bgamma^\circ)$
is asymptotically normally distributed. We accomplish this by analyzing
in turn each of the terms in the following decomposition:
\begin{eqnarray}
\lefteqn{{\bf 0}}&=&{\bf U}
(\hat{\bgamma},\hat{\Lambda}_0(\cdot,\hat{\bgamma}))\nonumber
\\ &=& {\bf U}(\bgamma^\circ,\Lambda_0^\circ) +
[{\bf U} (\bgamma^\circ,\hat{\Lambda}_0(\cdot,\bgamma^\circ))-{\bf
U} (\bgamma^\circ,\Lambda_0^\circ)] \nonumber\\ & & + [{\bf
U}(\hat{\bgamma},\hat{\Lambda}_0(\cdot,\hat{\bgamma}))- {\bf
U}(\bgamma^\circ,\hat{\Lambda}_0(\cdot,\bgamma^\circ))]. \nonumber
\end{eqnarray}
We show further that the covariance matrix of $\hat{\bgamma}$ can
be consistently estimated by the sandwich estimator
 \begin{equation}\label{eq:se}
  {\bf D}^{-1}(\hat{\bgamma})\{\hat{\bf V}(\hat{\bgamma})+\hat{\bf G}(\hat{\bgamma})+\hat{\bf C}(\hat{\bgamma})\}
  {\bf D}^{{-1}}(\hat{\bgamma})^T.
 \end{equation}
The matrix ${\bf D}$ consists of the derivatives of the $U_r$'s
with respect to the parameters $\bgamma$. ${\bf V}$ is the
asymptotic covariance matrix of ${\bf
U}(\bgamma^\circ,\Lambda_0^\circ)$, ${\bf G}$ is the asymptotic
covariance matrix of $[{\bf
U}(\bgamma^\circ,\hat{\Lambda}_0(\cdot,\bgamma^\circ))-{\bf
U}(\bgamma^\circ,\Lambda_0^\circ)]$,
 and ${\bf C}$ is the asymptotic covariance matrix between ${\bf U}(\bgamma^\circ,\Lambda_0^\circ)$
 and
$
[{\bf U}(\bgamma^\circ,\hat{\Lambda}_0(\cdot,\bgamma^\circ))-{\bf
U}(\bgamma^\circ,\Lambda_0^\circ)]$. The term ${\bf G}+{\bf C}$
reflects the added variance resulting from the need to estimate
the cumulative hazard function. All the above matrices are defined
explicitly in Section 4.4.

\section{Technical Conditions and Proofs}

This section presents the technical conditions we assume for the
asymptotic results and the proofs of these results.

\subsection{Technical Conditions}

In deriving the asymptotic properties of $\hat{\bgamma}$ we make
the following assumptions:

\begin{enumerate}
\item
The random vectors
$(T^0_{i1},\ldots,T^0_{im_i},C_{i1},\ldots,C_{im_i},{\bf
Z}_{i1},\ldots,{\bf Z}_{im_i}, W_i)$, $i=1,\ldots,n$, are
independent and identically distributed.
\item
There is a finite maximum follow-up time $\tau>0$, with
$\E[\sum_{j=1}^{m_i}Y_{ij}(\tau)]=y^*>0$ for all $i$.
\item
\begin{enumerate}
\item
Conditional on ${\bf Z}_{ij}$ and $W_i$, the censoring is independent and
noninformative of $W_i$ and $(\bbeta,\Lambda_0)$.
\item
$W_i$ is independent of ${\bf Z}_{ij}$ and of $m_i$.\
\end{enumerate}
\item
The frailty random variable $W_i$ has finite moments up to order
$(m+2)$, where $m$ is a fixed upper bound on $m_{i}$.
\item
${\bf Z}_{ij}$ is bounded.
\item
The parameter $\bgamma$ lies in a compact subset $\cG$ of $\real^{p+1}$
containing an open neighborhood of $\bgamma^\circ$.
\item
There exist $b>0$ and $C>0$ such that
$$
\lim_{w \rightarrow 0} w^{-(b-1)} f(w) = C.
$$
\item
The baseline hazard function $\lambda_0^\circ(t)$ is bounded over $[0,\tau]$
by some constant $\lambda_{max}$.
\item
The function $f^\prime(w;\theta) = (d/d\theta) f(w;\theta)$ is absolutely
integrable.
\item
The censoring distribution has at most finitely many jumps on $[0,\tau]$.
\item
 The matrix $[(\partial/\partial \bgamma)
\bU(\bgamma,\hat{\Lambda}_0(\cdot,\bgamma))]|_{\bgams=\bgams^\circ}$
is invertible with probability going to 1 as $n \rightarrow
\infty$.
\end{enumerate}

\subsection{Technical Preliminaries}

Since $\bbeta$ and ${\bf Z}_{ij}$ are bounded,
there exists a constant $\nu>0$ such that
\begin{equation}
\nu^{-1} \leq \exp(\bbeta^T {\bf Z}_{ij}) \leq \nu.
\label{eb}
\end{equation}
Now recall that $$ \psi_i( \bgamma,\Lambda,t) = \frac{\int
w^{N_i(t)+1} e^{-H_{i \cdot}(t)w} f(w) dw} {\int w^{N_i(t)}
e^{-H_{i \cdot}(t)w} f(w) dw}, $$ with $H_{i \cdot}(t)=H_{i
\cdot}(t,\bgamma,\Lambda) = \sum_{j=1}^{m_i} \Lambda(T_{ij} \wedge
t) \exp(\bbeta^T {\bf Z}_{ij})$ (here we define $H_{i \cdot}$ so
as to allow dependence on a general $\bgamma$ and $\Lambda$, which
will often not be explicitly indicated in the notation). Define
(for $0 \leq r \leq m$ and $h \geq 0$) $$ \psi^*(r,h) = \frac{\int
w^{r+1} e^{-hw} f(w) dw} {\int w^{r} e^{-hw} f(w) dw}. $$ Also
define $\psi_{min}^*(h) = \min_{0 \leq r \leq m} \psi^*(r,h)$ and
$\psi_{max}^*(h) = \max_{0 \leq r \leq m} \psi^*(r,h)$.
In the expression for $\psi^*(r,h)$, the numerator and denominator
are bounded above since $W$ is assumed to have finite
$(m+2)$-th moment. In addition, since $W$ is nondegenerate,
the numerator and denominator are strictly positive.
Thus $\psi_{max}^*(h)$
is finite and $\psi_{min}^*(h)$ is strictly positive.

\vspace*{0.5em}
\noindent {\bf Lemma 1:} The function $\psi^*(r,h)$ is decreasing in $h$.
Hence for all $\bgamma \in \cG$ and all $t$,
\begin{eqnarray}
\psi_i(\bgamma,\Lambda,t) & \leq & \psi_{max}^*(0), \label{huey} \\
\psi_i(\bgamma,\Lambda,t) & \geq & \psi_{min}^*(m \nu \Lambda(t)). \label{louie}
\end{eqnarray}
In addition, there exist $B>0$ and $\bar{h}>0$ such that,
for all $h \geq \bar{h}$,
\begin{equation}
\psi_{min}^*(h) \geq B h^{-1}.
\label{moe}
\end{equation}

\noindent {\bf Proof:}
We have
\begin{equation}
\frac{\partial}{\partial h} \psi^*(r,h)
= - \left[
\frac{\int w^{r+2} e^{-hw} f(w) dw}
{\int w^{r} e^{-hw} f(w) dw}
- \left (\frac{\int w^{r+1} e^{-hw} f(w) dw}
{\int w^{r} e^{-hw} f(w) dw} \right)^2 \right].
\label{psider}
\end{equation}
This is negative for all $h$, and so $\psi^*(r,h)$
is decreasing in $h$. Now $\psi_i(\bgamma,\Lambda,t) =
\psi^*(N_i(t),H_{i \cdot}(t))$. Since $0 \leq H_{i \cdot}(t) \leq m \nu \Lambda(t)$,
(\ref{huey}) and (\ref{louie}) follow.
As for (\ref{moe}), from a change of variable and Assumption 7,
$$
\lim_{h \rightarrow \infty} h \psi^*(r,h) =
\frac{\int_0^\infty v^{r+b} e^{-v} dv} {\int_0^\infty v^{r+b-1}
e^{-v} dv} = r + b.
$$
Now just take $\bar{h}$ large enough so that this
limit is obtained up to some factor, e.g.\ 1.01.

\vspace*{0.5cm}
\noindent {\bf Lemma 2:}
Define $\bar{\Lambda} = 1.03 e^{m\sigma} \bar{h}/(m \nu)$,
with $\sigma = 1.01 m \nu^2 /(By^*)$, with $\bar{h}$ and $B$ as above.
Then, with probability one, there exists $n^\prime$ such that, for all $t \in [0,\tau]$
and $\bgamma \in \cG$,
\begin{equation}
\hat{\Lambda}_0(t,\bgamma) \leq \bar{\Lambda} \quad \mbox{for } n \geq n^\prime,
\label{lamb}
\end{equation}
Thus, $\hat{\Lambda}_0(t,\bgamma)$
is naturally bounded, with no need to impose an upper bound artificially.

\noindent{\bf Proof:} To simplify the writing below, we will
suppress the argument $\bgamma$ in $\hat{\Lambda}_0(t,\bgamma)$.
Recall
$$
\Delta \hat{\Lambda}_0(\tau_k) =
\left[
 {\sum_{i=1}^{n}
 \psi_{i}(\bgamma,\hat{\Lambda}_0,\tau_{k-1})
 \sum_{j=1}^{m_i} Y_{ij}(\tau_k) \exp(\bbeta^T {\bf Z}_{ij})} \right]^{-1},
$$
where we now take $d_k=1$ since the survival time distribution is
assumed continuous. Using Lemma 1 and (\ref{eb}), we have
$$
\Delta \hat{\Lambda}_0(\tau_k) \leq
n^{-1} \nu
\psi_{min}^*(m \nu \hat{\Lambda}(\tau_{k-1}))
^{-1}
\left [ \frac{1}{n} \sum_{i=1}^n \sum_{j=1}^{m_i} Y_{ij}(\tau) \right]^{-1}.
$$
By the strong law of large numbers, there exists with probability one some
$n^*$ such that
\begin{equation}
\frac{1}{n} \sum_{i=1}^n \sum_{j=1}^{m_i} Y_{ij}(\tau)
\geq 0.999 y^*  \quad \mbox{for } n \geq n^*.
\label{yy}
\end{equation}
We thus have, for $n \geq n^*$,
\begin{equation}
\Delta \hat{\Lambda}_0(\tau_k) \leq
n^{-1} \left( \frac{1.01 \nu}{y^*} \right)
\psi_{min}^*(m \nu \hat{\Lambda}(\tau_{k-1}))^{-1}.
\label{del}
\end{equation}

Now, if $\hat{\Lambda}_0(t) \leq \bar{h}/(m\nu)$ for all $t$ then
we are done. Otherwise, there exists $k^\prime$ such that
$\hat{\Lambda}_0(\tau_k) \leq \bar{h}/(m\nu)$ for $k < k^\prime$
and $\hat{\Lambda}_0(\tau_k) \geq \bar{h}/(m\nu)$ for $k \geq
k^\prime$. Using the last inequality of Lemma 1, we obtain, for
$k>k^\prime$, $$ \Delta \hat{\Lambda}_0(\tau_k) \leq n^{-1} \sigma
\hat{\Lambda}_0(\tau_{k-1}), $$ or, in other words, $$
\hat{\Lambda}_0(\tau_k) \leq \left( 1 + \frac{\sigma}{n} \right)
\hat{\Lambda}_0(\tau_{k-1}). $$ Iterating the above inequality we
get $$ \hat{\Lambda}_0(\tau_{k^\prime+\ell}) \leq \left( 1 +
\frac{\sigma}{n} \right)^\ell \hat{\Lambda}_0(\tau_{k^\prime})
\leq \left( 1 + \frac{\sigma}{n} \right)^{mn}
\hat{\Lambda}_0(\tau_{k^\prime}) \leq 1.01 e^{m \sigma}
\hat{\Lambda}_0(\tau_{k^\prime}) $$ for $n$ large enough. But,
using (\ref{del}) and the fact that
$\hat{\Lambda}_0(\tau_{k^\prime-1}) \leq \bar{h}/(m\nu)$, we have
$$ \hat{\Lambda}_0(\tau_{k^\prime}) \leq \frac{\bar{h}}{m \nu} +
n^{-1} \left( \frac{1.01 \nu}{y^*} \right)
\psi_{min}^*(\bar{h})^{-1}, $$ which is less than $1.01
{\bar{h}}/({m \nu})$ for $n$ large enough. The desired conclusion
follows.

\vspace*{0.5cm}
\noindent {\bf Lemma 3:} We have $\sup_{s \in [0,\tau]}|\hat{\Lambda}_0(s,\bgamma^\circ)
-\hat{\Lambda}_0(s-,\bgamma^\circ)| \rightarrow 0 $ as $n \rightarrow \infty$,
as an immediate consequence of Lemma 2 and (\ref{del}).

\subsection{Consistency}

We now show the almost sure consistency of $\hat{\bbeta}$ and
$\hat{\Lambda}_0$. The argument is built on Claims~$\mbox{A-C}$ of
Section 3, which we prove below. Our argument follows
Zucker (2005, Appendix A.3).

\vspace*{0.5cm}
\noindent
{\bf Claim A:} $\hat{\Lambda}_0(t,\bgamma)$ converges a.s.\
to some function $\Lambda_0(t,\bgamma)$ uniformly in $t$ and $\bgamma$.

\noindent {\bf Proof:} Whenever a functional norm is
written below, the relevant uniform norm is intended.
We define $\Lambda_{max} = \max(\bar{\Lambda},\lambda_{max}\tau)$
and
$\psi^{**}(r,h) = \psi^*(r,h \wedge h_{max})$,
where $h_{max}$ = $m \nu \Lambda_{max}$. It is easy to see from (\ref{psider})
that $\psi^{**}(r,h)$ is Lipschitz continuous in $h$ (uniformly in~$r$). Recall
that $\psi_i(\bgamma,\Lambda,t) = \psi^*(N_i(t),H_{i \cdot}(t,\bgamma,\Lambda))$.
Lemma 2 implies that $H_{i \cdot}(t,\bgamma,\hat{\Lambda}_0(\cdot,\bgamma)) \leq h_{max}$
for all $t \in [0,\tau]$ and $\bgamma \in \cG$. Hence $\psi_i(\bgamma,\hat{\Lambda}_0(\cdot,\bgamma),t) =
\psi^{**}(N_i(t),H_{i \cdot}(t,\bgamma,\hat{\Lambda}_0(\cdot,\bgamma)))$.

Now define, for a general function $\Lambda$,
$$
\Xi_n(t,\bgamma,\Lambda)=\int_0^t
\frac{n^{-1} \sum_{i=1}^{n}\sum_{j=1}^{m_i}dN_{ij}(s)}
{n^{-1} \sum_{i=1}^{n}\sum_{j=1}^{m_i}\psi^{**}(N_i(s-),H_{i \cdot}(s-,\bgamma,\Lambda))
Y_{ij}(s)\exp(\bbeta^T {\bf Z}_{ij})}
$$
and
$$
\Xi(t,\bgamma,\Lambda)=\int_0^t
 \frac{\E [\sum_{j=1}^{m_i}\psi^*(N_i(s-),H_{i \cdot}(s-,\bgamma^\circ,\Lambda_0^\circ))
Y_{ij}(s)\exp(\bbeta^{\circ T}{\bf Z}_{ij})]}
{\E [\sum_{j=1}^{m_i}\psi^{**}(N_i(s-),H_{i \cdot}(s-,\bgamma,\Lambda))Y_{ij}(s-)
\exp(\bbeta^T {\bf Z}_{ij})]} \lambda_0^\circ(s) ds.
$$
By definition, $\hat{\Lambda}_0(t,\bgamma)$ satisfies the equation
\begin{equation}
\hat{\Lambda}_0(t,\bgamma) = \Xi_n(t,\bgamma,\hat{\Lambda}_0(\cdot,\bgamma)).
\label{lest}
\end{equation}

Next, define
$$
q_{\bgams}(s,\Lambda)
= \frac{\E [\sum_{j=1}^{m_i}\psi^*(N_i(s-),H_{i \cdot}(s-,\bgamma^\circ,\Lambda_0^\circ))
Y_{ij}(s)\exp(\bbeta^{\circ T}{\bf Z}_{ij})]}
{\E [\sum_{j=1}^{m_i}\psi^{**}(N_i(s-),H_{i \cdot}(s-,\bgamma,\Lambda))Y_{ij}(s)
\exp(\bbeta^T {\bf Z}_{ij})]} \lambda_0^\circ(s).
$$
This function is uniformly bounded by
$B^* = [\psi_{max}^*(0)/\psi_{min}^*(h_{max})]\lambda_{max}$.
Moreover,
by the Lipschitz continuity of $\psi^{**}(r,h)$ with respect to $h$,
it satisfies a Lipschitz-like condition of the form
$|q_{\bgams}(s,\Lambda_1)-q_{\bgams}(s,\Lambda_2)| \leq
K \sup_{0 \leq u \leq s} |\Lambda_1(u)-\Lambda_2(u)|$.
Hence, by mimicking the argument of Hartman (1973, Theorem 1.1), we find that the
equation $\Lambda(t) = \Xi(t,\bgamma,\Lambda)$ has a unique solution, which we denote
by $\Lambda_0(t,\bgamma)$. The claim then is that $\hat{\Lambda}_0(t,\bgamma)$
converges almost surely (uniformly in $t$ and $\bgamma$) to $\Lambda_0(t,\bgamma)$.
Though it may be possible to prove this claim directly, we shall use a convenient indirect argument.

Define $\tilde{\Lambda}_0^{(n)}(t,\bgamma)$ to be a modified version of
$\hat{\Lambda}_0(t,\bgamma)$ defined by linear interpolation between
the jumps.
Lemma 3 implies that, with probability one,
\begin{equation}
\sup_{t,\bgams} |\tilde{\Lambda}_0^{(n)}(t,\bgamma)
- \hat{\Lambda}_0(t,\bgamma)| \rightarrow 0, \label{Lone}
\end{equation}
and thus
\begin{equation}
\sup_{t,\bgams} |\Xi_n(t,\bgamma,\tilde{\Lambda}_0(t,\bgamma))-
\Xi_n(t,\bgamma,\hat{\Lambda}_0(t,\bgamma))| \rightarrow 0. \label{Ltwo}
\end{equation}
Lemma 2 shows that the family $\cL = \{ \tilde{\Lambda}_0^{(n)}(t,\bgamma)$,
$n \geq n^\prime \}$ is uniformly bounded. We can show further that $\cL$ is
equicontinuous. This is done as follows.

Recall that $N_i(t)=\sum_{j=1}^{m_i} N_{ij}(t)$.
Write $\bar{N}(t) = n^{-1} \sum_{i=1}^n \sum_{j=1}^{m_i} N_{ij}(t)$.
We have $\bar{N}(t) \rightarrow \E[N_i(t)]$ as $n \rightarrow \infty$
uniformly in $t$ with probability one, with
$$
\E[N_i(t)] =
\int_0^t \E \left[ \sum_{j=1}^{m_i}\psi^*(N_i(s-),H_{i \cdot}(s-,\bgamma^\circ,\Lambda_0^\circ))
Y_{ij}(s)\exp(\bbeta^{\circ T}{\bf Z}_{ij}) \right] \lambda_0^\circ(s) ds.
$$
In view of this and (\ref{yy}) there exists a probability-one set of realizations
$\Omega^*$ on which the following holds: for any given $\epsilon>0$, we can find
$n^{\prime \prime}(\epsilon)$ such that $\sup_t |\bar{N}(t)-E[N_i(t)]| \leq \epsilon/(4B^\circ)$
for all $n \geq n^{\prime \prime}(\epsilon)$, where $B^\circ = 1.01 \nu/
[\psi_{min}^*(h_{max})y^*]$.
In consequence, for all $t$ and $u$ with $u<t$, we find that
$$
 \hat{\Lambda}_0(t,\bgamma) - \hat{\Lambda}_0(u,\bgamma)
 = \int_u^t
\frac{n^{-1} \sum_{i=1}^{n}\sum_{j=1}^{m_i}dN_{ij}(s)}
{n^{-1} \sum_{i=1}^{n}\sum_{j=1}^{m_i}\psi^{**}(N_i(s-),H_{i \cdot}(s-,\bgamma,\Lambda))
Y_{ij}(s)\exp(\bbeta^T {\bf Z}_{ij})}
$$
satisfies
\begin{equation}
\hat{\Lambda}_0(t,\bgamma) - \hat{\Lambda}_0(u,\bgamma)
\leq B^* (t-u) + \frac{\epsilon}{2} \quad \mbox{for all } n \geq n^{\prime \prime}(\epsilon).
\label{equi}
\end{equation}
Moreover, it is easy to see that $\hat{\Lambda}_0(t,\bgamma)$ is
Lipschitz continuous in $\bgamma$ with Lipschitz constant
$C^*$, say, that is independent of $t$.

These two results imply that $\cL$ is equicontinuous. This is seen as follows.
For given $\epsilon$, we need to find $\delta_1^*$ and $\delta_2^*$ such that
$|\tilde{\Lambda}_0^{(n)}(t,\bgamma) -\tilde{\Lambda}_0^{(n)}(u,\bgamma)| \leq \epsilon$ whenever $|t-u|
\leq \delta_1^*$ and $|\tilde{\Lambda}_0^{(n)}(t,\bgamma) -\tilde{\Lambda}_0^{(n)}(t,\bgamma^\prime)| \leq
\epsilon$ whenever $\|\bgamma-\bgamma^\prime \| \leq \delta_2^*$. The latter is easily obtained using the
Lipschitz continuity of $\hat{\Lambda}_0(t,\bgamma)$ with respect to $\bgamma$. As for the former,
for $n \geq n^{\prime \prime}(\epsilon)$ this can be accomplished using (\ref{equi}), while for $n$ in
the finite set $n^\prime \leq n < n^{\prime \prime}(\epsilon)$ this can be accomplished using the
fact that the function $\tilde{\Lambda}_0^{(n)}(t,\bgamma)$ is uniformly continuous on $[0,\tau]$ for
every given $n$.

We have thus shown that $\cL$ is (almost surely)
a relatively compact set in the space $C([0,\tau] \times \cG)$.

Next, define
\begin{eqnarray*}
A(\bgamma,\Lambda,s) & = &
\frac{1}{n} \sum_{i=1}^{n}\sum_{j=1}^{m_i}\psi^{**}(N_i(s-),H_{i \cdot}(s-,\bgamma,\Lambda))
Y_{ij}(s)\exp(\bbeta^T {\bf Z}_{ij}), \\
a(\bgamma,\Lambda,s) & = &
\E \left[\sum_{j=1}^{m_i}\psi^{**}(N_i(s-),H_{i \cdot}(s-,\bgamma,\Lambda))
Y_{ij}(s)\exp(\bbeta^T {\bf Z}_{ij}) \right].
\end{eqnarray*}
For any fixed continuous $\Lambda$, the functional strong law of large numbers of
Andersen \& Gill (1982, Appendix III) implies that
\begin{equation}
\sup_{s,\bgams} |A(\bgamma,\Lambda,s) - a(\bgamma,\Lambda,s)| \rightarrow 0
\hspace*{1em} \mbox{a.s.}
\label{ac}
\end{equation}
Here we need the following more complex result:
\begin{equation}
\sup_{s,\bgams} |A(\bgamma,\tilde{\Lambda}^{(n)},s) - a(\bgamma,\tilde{\Lambda}^{(n)},s)| \rightarrow 0
\hspace*{1em} \mbox{a.s.}
\label{acn}
\end{equation}
The proof of (\ref{acn}) is lengthy; we give the details in Section 4.5 below.
In outline form, the proof involves two steps: (1) showing that, for any given $\epsilon > 0$,
we can define an appropriate \it finite \rm class $\cL_\epsilon^*$ of functions $\Lambda$
such that $\tilde{\Lambda}^{(n)}$ can be suitably approximated by some member of the class;
(2) applying the result (\ref{ac}), which will hold uniformly over the finite class.

Given (\ref{acn}) and the a.s. uniform convergence of
$\bar{N}(t)$ to $\E[N_i(t)]$, we can infer that
\begin{equation}
\sup_{t,\bgams}
|\Xi_n(t,\bgamma,\tilde{\Lambda}_0^{(n)}(t,\bgamma))
- \Xi(t,\bgamma,\tilde{\Lambda}_0^{(n)}(t,\bgamma))| \rightarrow 0
\hspace*{1em} \mbox{a.s.}
\label{supX}
\end{equation}
The result (\ref{supX}) is easily obtained by adapting the
argument of Aalen (1976, Lemma~6.1), using the
equicontinuity of $\cL$. It is here that we use
Assumption~10, for the adaptation of Aalen's argument requires
$a(\bgamma,\Lambda,s)$ to be piecewise continuous with finite left
and right limits at each point of discontinuity.

From (\ref{lest}), (\ref{Lone}), (\ref{Ltwo}), and (\ref{supX}) it
follows that any limit point of $\{
\tilde{\Lambda}_0^{(n)}(t,\bgamma) \}$ must satisfy the equation
$\Lambda=\Xi(t,\bgamma,\Lambda)$. Since $\Lambda_0(t,\bgamma)$ is
the unique solution of this equation, it is the unique limit point
of $\{ \tilde{\Lambda}_0^{(n)}(t,\bgamma) \}$. Thus $\{
\tilde{\Lambda}_0^{(n)}(t,\bgamma) \}$ is a sequence in a compact
set with unique limit point $\Lambda_0(t,\bgamma)$. Hence
$\tilde{\Lambda}_0^{(n)}(t,\bgamma)$ converges a.s.\ uniformly in
$t$ and $\bgamma$ to $\Lambda_0(t,\bgamma)$. In view of
(\ref{Lone}), the same holds of $\hat{\Lambda}_0(t,\bgamma)$,
which is the desired result.
Note that $\Lambda_0(\cdot,\bgamma^\circ)=\Lambda_0^\circ(\cdot)$
since $\Lambda_0^\circ$ trivially solves the equation $\Lambda=\Xi(t,\bgamma^\circ,\Lambda)$.
\\

\noindent
{\bf Claim B:} With $\bu(\bgamma,\Lambda_0(\cdot,\bgamma))
= \E[ \bU(\bgamma,\Lambda_0(\cdot,\bgamma)) ]$, we have
$\bU(\bgamma,\hat{\Lambda}_0(\cdot,\bgamma)) \rightarrow
\bu(\bgamma,\Lambda_0(\cdot,\bgamma))$ uniformly in
$\bgamma \in \cG$ with probability one.

\noindent {\bf Proof:}
Since $\bU(\bgamma,\Lambda_0(\cdot,\bgamma))$ is the mean of iid terms,
the functional strong law of numbers of Andersen \& Gill (1982, Appendix III) implies that
$\bU(\bgamma,\Lambda_0(\cdot,\bgamma))$ converges uniformly in $\bgamma$ almost surely
to $\bu(\bgamma,\Lambda_0(\cdot,\bgamma))$. It remains only to show that
\begin{equation}
\sup_{\bgams} \|\bU(\bgamma,\hat{\Lambda}_0(\cdot,\bgamma)) - \bU(\bgamma,\Lambda_0(\cdot,\bgamma))\|
\rightarrow \bzro
\label{fuzz}
\end{equation}
almost surely. The structure of $\bU(\bgamma,\Lambda)$ reveals that there
exists some constant $C^\circ$ (independent of $\bgamma$) such that
$\|\bU(\bgamma,\Lambda_1)-\bU(\bgamma,\Lambda_2)\| \leq C^\circ \| \Lambda_1 - \Lambda_2 \|.$
>From this along with Claim A, (\ref{fuzz}) follows.
\\

\noindent {\bf Claim C:} There exists a unique consistent root to
${\bf U}(\hat{\bgamma},\hat{\Lambda}_0(\cdot,\hat{\bgamma}))=\bzro$.

\noindent {\bf Proof:}
We apply Foutz's (1977) consistency theorem for maximum likelihood type
estimators. The following conditions must be established:
\begin{flushleft}
{\bf F1.}
$\partial \bU(\bgamma,\hat{\Lambda}_0(\cdot,\bgamma))/\partial
\bgamma$ exists and is continuous in an open neighborhood about
$\bgamma^\circ$. \\
{\bf F2.}
The convergence of $\partial
\bU(\bgamma,\hat{\Lambda}_0(\cdot,\bgamma)) /\partial \bgamma$ to
its limit is uniform in open neighborhood of $\bgamma^\circ$. \\
{\bf F3.}
$\bU(\bgamma^\circ,\hat{\Lambda}_0(\cdot,\bgamma^\circ)) \rightarrow \bzro$ as
$n \rightarrow \infty$. \\
{\bf F4.}
The matrix $-[\partial
\bU(\bgamma,\hat{\Lambda}_0(\cdot,\bgamma))/\partial
\bgamma]|_{\bgams=\bgams^\circ}$ is invertible with probability going to 1
as $n \rightarrow \infty$. (In Foutz's paper, the matrix in question is symmetric,
and so he stated the condition in terms of positive definiteness. But his proof,
which is based on the inverse function theorem, shows that the basic condition needed
is invertibility.)
\end{flushleft}
It is easily seen that Condition F1 holds. Given Assumptions 2, 4,
and 5, Condition F2 follows from the previously-cited functional
law of large numbers. As for Condition F3, in Claim B we showed
that $\bU(\bgamma,\Lambda_0(\cdot,\bgamma))$ converges a.s.\
uniformly to $\bu(\bgamma,\Lambda_0(\cdot,\bgamma)) = \E
[\bU(\bgamma,\Lambda_0(\cdot,\bgamma))]$. We noted already that
$\Lambda_0(\cdot,\bgamma^\circ)=\Lambda_0(\cdot)$. Thus we
need only show that $\E [\bU(\bgamma^\circ,\Lambda_0)]=~{\bf 0}$.
Since $\bU$ is a score function derived from a classical iid
likelihood, this result follows from classical likelihood theory.
Condition F4 has been assumed in Assumption 11.
With Conditions F1-F4 established, the result follows.

\subsection{Asymptotic Normality}
To show that $\hat{\bgamma}$ is asymptotically normally
distributed, we write
\begin{eqnarray}
\lefteqn{{\bf 0}}&=&\bU(\hat{\bgamma},\hat{\Lambda}_0(\cdot,\hat{\bgamma}))\nonumber
\\ &=& \bU(\bgamma^\circ,\Lambda_0^\circ) +
[\bU(\bgamma^\circ,\hat{\Lambda}_0(\cdot,\bgamma^\circ))-\bU(\bgamma^\circ,\Lambda_0^\circ)]
\nonumber\\ & & + \,
[\bU(\hat{\bgamma},\hat{\Lambda}_0(\cdot,\hat{\bgamma}))-
\bU(\bgamma^\circ,\hat{\Lambda}_0(\cdot,\bgamma^\circ))] \nonumber
\end{eqnarray}
In the following we consider each of the above terms of the
right-hand side of the equation.


\bsh
\underline{Step I}
\esh

We can write $\bU(\bgamma^\circ, \Lambda_0^\circ)
= n^{-1} \sum_{i=1}^n \bxi_i$,
where $\bxi_i$ is a $(p+1)$-vector with $r$-th element, $r=1,\ldots,p$,
given by
$$
\xi_{ir}  =  \sum_{j=1}^{m_i} \delta_{ij} Z_{ijr} -
 \frac{\left[\sum_{j=1}^{m_i}{H}_{ij}(\tau)
 Z_{ijr}\right] \int w^{N_{i.}(\tau)+1} \exp\{-w \{H_{i.}(\tau)\}f(w;{\theta})dw}
 {\int w^{N_{i.}(\tau)} \exp\{-w {H}_{i.}(\tau)\}f(w;{\theta})dw}
$$
and $(p+1)$-th element given by
$$
\xi_{i(p+1)} =  \frac{\int w^{N_{i.}(\tau)}
\exp\{-w{H}_{i.}(\tau)\} f'(w;{\theta})dw}
 {\int w^{N_{i.}(\tau)} \exp\{-w {H}_{i.}(\tau)\}f(w;{\theta})dw}.
$$
Thus $\bU(\bgamma^\circ, \Lambda_0^\circ)$ is the
mean of the iid mean-zero random vectors $\bxi_i$. It hence follows
from the central limit theorem that $n^{\half} \bU(\bgamma^\circ, \Lambda_0^\circ)$
is asymptotically mean-zero multivariate normal. To estimate the covariance
matrix, let $\bxi_i^*$ be the counterpart of $\bxi_i$ with estimates of $\bgamma$
and $\Lambda_0$ substituted for the true values. Then an empirical
estimator of the covariance matrix is given by
 $\hat{{\bf V}}(\hat{\bgamma}) = n^{-1} \sum_{i=1}^{n}
 \bxi_i^{*}\bxi_i^{*T}$.
This is a consistent estimator of the covariance matrix
since $\hat{\Lambda}_0(t,\bgamma)$ converges to $\Lambda_0(t,\bgamma)$ a.s.\
uniformly in $t$ and $\bgamma$  (Claim A), and $\hat{\bgamma}$ is a
consistent estimator of $\bgamma^\circ$ (Claim C).


\bsh
\underline{Step II}
\esh

Let $\hat{U}_r=U_r(\bgamma^\circ,\hat{\Lambda}_0)$, $r=1,\ldots,p$,
and $\hat{U}_{p+1}=U_{p+1}(\bgamma^\circ,\hat{\Lambda}_0)$ (in this
segment of the proof, when we write $(\bgamma^\circ,\hat{\Lambda}_0)$
the intent is to signify $(\bgamma^\circ,\hat{\Lambda}_0(\cdot,\bgamma^\circ))$.
First order Taylor expansion of $\hat{U}_r$ about $\Lambda_0^\circ$,
$r=1,\ldots,p+1$, gives
$$
\hspace*{-6cm}
n^{1/2}\{U_r(\bgamma^\circ,\hat{\Lambda}_0)-U_r(\bgamma^\circ,{\Lambda}_0^\circ)\}
$$
\begin{equation}\label{eq:tayloru}
\hspace*{1cm} = n^{-1/2} \sum_{i=1}^{n} \sum_{j=1}^{m_i}
Q_{ijr}(\bgamma^\circ,\Lambda^\circ,T_{ij})
\{\hat{\Lambda}_0(T_{ij},\bgamma^\circ)-\Lambda_0^\circ(T_{ij})\} +
o_p(1),
\end{equation}
where
\begin{eqnarray*}
Q_{ijr}(\bgamma^\circ,\Lambda^\circ,T_{ij})&=& -\left\{
 \frac{\phi_{2i}(\bgamma^\circ,\Lambda_0^\circ,\tau)}{\phi_{1i}(\bgamma^\circ,\Lambda_0^\circ,\tau)}
 R_{ij}^* Z_{ijr}
 -
 \frac{\phi_{3i}(\bgamma^\circ,\Lambda_0^\circ,\tau)}{\phi_{1i}(\bgamma^\circ,\Lambda_0^\circ,\tau)}
 R_{ij}^*\sum_{j=1}^{m_i} H_{ij}(T_{ij}) Z_{ijr}
 \right. \nonumber \\
& & \left.
\hspace*{1cm} + \, \frac{\phi_{2i}^2(\bgamma^\circ,\Lambda_0^\circ,\tau)}
 {\phi_{1i}^2(\bgamma^\circ,\Lambda_0^\circ,\tau)}
 R_{ij}^* \sum_{j=1}^{m_i} H_{ij}(T_{ij})Z_{ijr}  \right\}
\end{eqnarray*}
for $r=1,\ldots,p$, and
\begin{eqnarray*}
Q_{ij(p+1)}(\bgamma^\circ,\Lambda^\circ,T_{ij}) =  R_{ij}^* \left\{
\frac{\phi_{2i}(\bgamma^\circ,\Lambda_0^\circ,\tau)
 \phi_{1i}^{(\theta)}(\bgamma^\circ,\Lambda_0^\circ,\tau)}
 {\phi_{1i}^2(\bgamma^\circ,\Lambda_0^\circ,\tau)}
 -\frac{ \phi_{2i}^{(\theta)}(\bgamma^\circ,\Lambda_0^\circ,\tau) }
 {\phi_{1i}(\bgamma^\circ,\Lambda_0^\circ,\tau)}
\right\},
\end{eqnarray*}
with $R_{ij}^*=\exp(\bbeta^T {\bf Z}_{ij})$ and
 \begin{eqnarray*}
 \phi_{ki}^{(\theta)}(\bgamma,\Lambda_0,t)
 =
 \int w^{N_{i.}(t)+(k-1)} \exp\{-w H_{i.}(t)\} f'(w) dw,
\quad k=1,2.
 \end{eqnarray*}
The validity of the approximation (\ref{eq:tayloru}) can be seen by an argument
similar to that used in connection with (\ref{eq:yapp}) below.

Given the intensity process (\ref{eq:inten}), the process
 \begin{eqnarray*}
 M_{ij}(t)=N_{ij}(t)-\int_0^t
 \lambda_0(u) \exp(\bbeta^{\circ T} {\bf Z}_{ij}) Y_{ij}(u) \psi_{i}(\bgamma^\circ,\Lambda_0^\circ,u-) du
 \end{eqnarray*}
is a mean zero martingale with respect to the filtration
${\mathcal F}_t$. Also, by Lemma 3, we have that
$\sup_{s \in [0,\tau]} |\hat{\Lambda}_0(s,\bgamma^\circ)-\hat{\Lambda}_0(s-,\bgamma^\circ)|$
converges to zero. Thus, replacing $s-$ by $s$ we obtain the following
approximation, uniformly over $t \in [0,\tau]$:
 \begin{eqnarray}\label{eq:lamapp}
\hat{\Lambda}_0(t,\bgamma^\circ)-\Lambda_0^\circ(t) & \approx & \frac{1}{n}
\int_0^t \{\mathcal{Y}(s,\Lambda_0^\circ)\}^{-1} \sum_{i=1}^{n}
 \sum_{j=1}^{m_i} dM_{ij}(s) \nonumber\\
 &+& \frac{1}{n} \int_0^t \left[\{\mathcal{Y}(s,\hat{\Lambda}_0)\}^{-1}
 - \{\mathcal{Y}(s,\Lambda_0^\circ)\}^{-1}\right]\sum_{i=1}^{n}
 \sum_{j=1}^{m_i} dN_{ij}(s),
 \end{eqnarray}
where
 $$
 \mathcal{Y}(s,\Lambda) = \frac{1}{n}  \sum_{i=1}^{n}
 \psi_i(\bgamma^\circ,\Lambda,s)\sum_{j=1}^{m_i} Y_{ij}(s)
 \exp(\bbeta^{\circ T} {\bf Z}_{ij}).
 $$

Now let
${\mathcal W}(s,r)=\{{\mathcal Y}(s,\Lambda_0^\circ+r\Delta)\}^{-1}$
with $\Delta=\hat{\Lambda}_0-\Lambda_0^\circ$. Define
$\dot{\mathcal W}$ and $\ddot{\mathcal W}$ as the first and second
derivative of ${\mathcal W}$ with respect to $r$, respectively.
Then, computing the necessary derivatives and carrying out
a first order Taylor expansion of ${\mathcal W}(s,r)$
around $r=0$ evaluated at $r=1$ with Lagrange remainder
(Abramowitz \& Stegun, 1972, p. 880), we get
$$
\hspace*{-38em} \lefteqn{\{{\mathcal Y}(s,\hat{\Lambda}_0)\}^{-1} - \{{\mathcal
Y}(s,{\Lambda}_0^\circ)\}^{-1}
 =\dot{\mathcal
W}(s,0)+\frac{1}{2}\ddot{\mathcal W}(s,\tilde{r}(s))}
$$
\begin{equation}
= -\frac{1}{n} \sum_{i=1}^{n}\sum_{j=1}^{m_i} \left[
 \frac{R_{i.}(s)\eta_{1i}(0,s)}{\{{\mathcal Y}(s,{\Lambda}_0^\circ)\}^2}
 -\frac{1}{2}h_i(\tilde{r}(s),s) \right]
  \exp(\bbeta^T {\bf Z}_{ij})
  \{ \hat{\Lambda}_0(T_{ij}\wedge s) - {\Lambda}_0^\circ(T_{ij}\wedge s)\},
  \label{eq:yapp}
\end{equation}
where $R_{ij}(u)=\exp(\bbeta^T {\bf Z}_{ij}) Y_{ij}(u)$,
$R_{i.}(u)=\sum_{j=1}^{m_i} R_{ij}(u)$, $\tilde{r}(s) \in [0,1]$,
\begin{eqnarray*}
\eta_{1i}(r,s)=\frac{\phi_{3i}(\bgamma^\circ,{\Lambda}_0^\circ+r\Delta,s)}
{\phi_{1i}(\bgamma^\circ,{\Lambda}_0^\circ+r\Delta,s)}-\left\{
\frac{\phi_{2i}(\bgamma^\circ,{\Lambda}_0^\circ+r\Delta,s)}{\phi_{1i}(\bgamma^\circ,{\Lambda}_0^\circ+r\Delta,s)}
\right\}^2,
\end{eqnarray*}
and $h_i(r,s)$ is as defined in Section 4.6 below, and shown there
to be $o(1)$ uniformly in $r$ and $s$.

Let $\eta_{1i}(s)=\eta_{1i}(0,s)$. Plugging (\ref{eq:yapp}) into
(\ref{eq:lamapp}) we get
 \begin{eqnarray*}
 \lefteqn{\hat{\Lambda}_0(t,\bgamma^\circ)-\Lambda_0^\circ(t)  \approx
 n^{-1}\int_0^t \{\mathcal{Y}(s,\Lambda_0^\circ)\}^{-1} \sum_{i=1}^{n}
 \sum_{j=1}^{m_i} dM_{ij}(s)} \\
 && -n^{-2}\int_0^t
 \sum_{k=1}^{n}\sum_{l=1}^{m_k}
 \frac{I(T_{kl}>s) R_{k.}(s) \eta_{1k}(s)}
 {\{\mathcal{Y}(s,\Lambda_0^\circ)\}^{2}}
 \exp(\bbeta^T {\bf Z}_{kl})
 \{\hat{\Lambda}_0(s)-\Lambda_0^\circ(s)\}\sum_{i=1}^{n}
 \sum_{j=1}^{m_i} dN_{ij}(s)  \\
 && - n^{-2}\int_0^t
\sum_{k=1}^{n}\sum_{l=1}^{m_k}\frac{I(T_{kl} \leq s) R_{k.}(s)
 \eta_{1k}(s)}
 {\{\mathcal{Y}(s,\Lambda_0^\circ)\}^{2}}
 \exp(\bbeta^T {\bf Z}_{kl})
 \{\hat{\Lambda}_0(T_{kl})-\Lambda_0^\circ(T_{kl})\}\sum_{i=1}^{n}
 \sum_{j=1}^{m_i} dN_{ij}(s) \\
 && + \, n^{-2}\int_0^t
 \sum_{k=1}^{n}\sum_{l=1}^{m_k} \frac{1}{2} h_{k}(\tilde{r}(s),s) \exp(\bbeta^T {\bf Z}_{kl})
 \{\hat{\Lambda}_0(T_{kl})-\Lambda_0^\circ(T_{kl})\}\sum_{i=1}^{n}
 \sum_{j=1}^{m_i} dN_{ij}(s).
\end{eqnarray*}
The third term of the above equation can be written, by interchanging the
order of integration, as
$$
n^{-2}\sum_{k=1}^{n}\sum_{l=1}^{m_k}\sum_{i=1}^{n}\sum_{j=1}^{m_i}
\int_0^t \frac{R_{k.}(s) \eta_{1k}(s)}
{\{\mathcal{Y}(s,\Lambda_0^\circ)\}^{2}} \exp(\bbeta^T {\bf Z}_{kl})
 \left[ \int_0^s
\{\hat{\Lambda}_0(u)-\Lambda_0^\circ(u) \} d\tilde{N}_{kl}(u)\}
\right] dN_{ij}(s)
$$
$$
 = \int_0^{t}
  \{\hat{\Lambda}_0(s)-\Lambda_0^\circ(s)\}
\sum_{i=1}^{n}\sum_{j=1}^{m_i} \Omega_{ij}(s,t) d\tilde{N}_{ij}(s),
$$
where
$\tilde{N}_{ij}(t)=I(T_{ij} \leq t)$ and
 $$
 \Omega_{ij}(s,t)=n^{-2}\int_s^t \{\mathcal{Y}(u,\Lambda_0^\circ)\}^{-2}
 R_{i.}(u) \eta_{1i}(u)\exp(\bbeta^T {\bf Z}_{ij}) \sum_{k=1}^n \sum_{l=1}^{m_k}
 dN_{kl}(u).
 $$
Hence we get
 \begin{eqnarray}\label{eq:lamapp2}
\hat{\Lambda}_0(t,\bgamma^\circ)-\Lambda_0^\circ(t)  &\approx&
 n^{-1}\int_0^t \{\mathcal{Y}(s,\Lambda_0^\circ)\}^{-1} \sum_{i=1}^{n}
 \sum_{j=1}^{m_i} dM_{ij}(s) \nonumber\\
 &-&  \int_0^t
 \{\hat{\Lambda}_0(s,\bgamma^\circ)-\Lambda_0^\circ(s)\}\sum_{i=1}^{n}
 \sum_{j=1}^{m_i}\{\delta_{ij} \Upsilon(s)+\Omega_{ij}(s,t) +o(n^{-1})\}
 d\tilde{N}_{ij}(s) \nonumber
 \end{eqnarray}
where
 $$
\Upsilon(s) =n^{-2}\{\mathcal{Y}(s,\Lambda_0^\circ)\}^{-2}
\sum_{k=1}^{n} \sum_{l=1}^{m_k} I(T_{kl}>s) R_{k.}(s)
\eta_{1k}(s)\exp(\bbeta^T {\bf Z}_{kl}).
 $$
The $o(n^{-1})$ is uniform in $t$ (see Sec.~4.6 below) and will be
dominated by $\Omega$ and $\Upsilon$, which are of order $n^{-1}$.
Hence the $o(n^{-1})$ term can be ignored.

An argument similar to that
of Yang \& Prentice (1999) and Zucker (2005) now yields
the martingale representation
\begin{eqnarray}\label{eq:martin}
\hat{\Lambda}_0(t,\bgamma^\circ)-\Lambda_0^\circ(t)
 &\approx&
 \frac{1}{n \hat{p}(t)} \int_0^t
 \frac{\hat{p}(s-) \sum_{i=1}^{n}\sum_{j=1}^{m_i} dM_{ij}(s)}
 {\mathcal{Y}(s,\Lambda_0^\circ)},
\end{eqnarray}
where
\begin{eqnarray*}
\hat{p}(t)=\prod_{s \leq t} \left[
 1+\sum_{i=1}^{n} \sum_{j=1}^{m_i}
 \{\delta_{ij} \Upsilon(s)+\Omega_{ij}(s,t) \} d\tilde{N}_{ij}(s)
 \right].
\end{eqnarray*}
Based on (\ref{eq:tayloru}), we can write
 \begin{eqnarray*}
U_r(\bgamma^\circ,\hat{\Lambda}_0)-U_r(\bgamma^\circ,{\Lambda}_0^\circ)
\approx n^{-1} \sum_{i=1}^{n} \sum_{j=1}^{m_i}
 \int_0^{\tau} Q_{ijr}(\bgamma^\circ,\Lambda_0^\circ,s)
\{\hat{\Lambda}_0(s,\bgamma^\circ)-\Lambda_0^\circ(s)\} d
\tilde{N}_{ij}(s).
\end{eqnarray*}
Plugging the martingale representation (\ref{eq:martin}) into the
above equation and carrying out some more algebra (again involving
an interchange of integrals) gives
\begin{eqnarray}
U_r(\bgamma^\circ,\hat{\Lambda}_0)-U_r(\bgamma^\circ,{\Lambda}_0^\circ)
& & \nonumber \\
& \hspace*{-6cm}  \approx & \hspace*{-3cm}
n^{-1}
\int_0^{\tau} \pi_r(s,\bgamma^\circ,\Lambda_0^\circ) \frac{\hat{p}(s-) \sum_{k=1}^{n}\sum_{l=1}^{m_k} dM_{kl}(s)}
{\mathcal{Y}(s,\Lambda_0^\circ)},
\label{eq:mart}
\end{eqnarray}
where
 $$
\pi_r(s,\bgamma,\Lambda_0) = n^{-1} \int_s^{\tau}
 \frac{\sum_{i=1}^{n} \sum_{j=1}^{m_i}
 Q_{ijr}(\bgamma,\Lambda_0,t)
 d\tilde{N}_{ij}(t)}{\hat{p}(t)}.
 $$
Therefore, $n^{1/2}[{\bf
U}(\bgamma^\circ,\hat{\Lambda}_0(\cdot,\bgamma^\circ))-{\bf
U}(\bgamma^\circ,\Lambda_0^\circ(\cdot,\bgamma^\circ))]$
 is asymptotically mean zero multivariate normal with covariance
 matrix that can be consistently estimated by
 \begin{eqnarray*}
 G_{rl}(\hat{\bgamma}) = n^{-1} \int_0^{\tau} \pi_r(s,\hat{\bgamma},\hat{\Lambda}_0)
  \pi_l(s,\hat{\bgamma},\hat{\Lambda}_0) \{\hat{p}(s-)\}^2
 \frac{\sum_{i=1}^{n} \sum_{j=1}^{m_i} dN_{ij}(s)}
 {\{\mathcal{Y}(s,\hat{\Lambda}_0)\}^{2}
 }
 \end{eqnarray*}
 for $r,l=1,\ldots,p+1$.

\bsh
\underline{Step III}
\esh

We now examine the sum of
$\bU(\bgamma^\circ, \Lambda_0^\circ)$ and
$\bU(\bgamma^\circ,\hat{\Lambda}_0(\cdot,\bgamma^\circ))-\bU(\bgamma^\circ,
 \Lambda_0^\circ)$. From (\ref{eq:mart}), we have
 \begin{eqnarray*}
U_r(\bgamma^\circ,\hat{\Lambda}_0(\cdot,\bgamma^\circ))-U_r(\bgamma^\circ,\Lambda_0^\circ)
\approx n^{-1}\int_0^{\tau} \alpha_r(s) \sum_{k=1}^n
\sum_{l=1}^{m_k} dM_{kl}(s) =\frac{1}{n}\sum_{k=1}^n \mu_{kr},
 \end{eqnarray*}
where $\alpha_r(s)$ is the limiting value of
$\pi_r(s,\bgamma^\circ,\Lambda_0^\circ) \hat{p}(s-)
/\mathcal{Y}(s,{\Lambda}_0^\circ)$
and $\mu_{kr}$ is defined as
\begin{eqnarray*}
\mu_{kr}=\int_0^{\tau} \alpha_r(s) \sum_{l=1}^{m_k} dM_{kl}(s).
\end{eqnarray*}
Arguments in Yang and Prentice (1999, Appendix A) can
be used to show that $\hat{p}(s-)$ has a limit. Also, clearly
$\E[\mu_{kr}]=0$.

We thus have
\begin{eqnarray*}
U_r(\bgamma^\circ,
\Lambda_0^\circ)+[U_r(\bgamma^\circ,\hat{\Lambda}_0(\cdot,\bgamma^\circ))-U_r(\bgamma^\circ,
 \Lambda_0^\circ)]
 \approx
\frac{1}{n} \sum_{i=1}^n (\xi_{ir} + \mu_{ir}),
\end{eqnarray*}
which is a mean of $n$ iid random variables. Hence
$n^{1/2}\{U_r(\bgamma^\circ,
\Lambda_0^\circ)+[U_r(\bgamma^\circ,\hat{\Lambda}_0(\cdot,\bgamma^\circ))-U_r(\bgamma^\circ,
 \Lambda_0^\circ)]\}$
is asymptotically normally distributed. The covariance matrix may
be estimated by
$\hat{\bV}(\hat{\bgamma})+\hat{\bG}(\hat{\bgamma})+\hat{\bC}(\hat{\bgamma})$,
where
\begin{eqnarray*}
\hat{C}_{rl}(\hat{\bgamma})=
\frac{1}{n}
\sum_{i=1}^{n} (\xi_{ir}^{*}
\mu_{il}^{*} + \xi_{il}^{*} \mu_{ir}^{*}), \quad
r,l=1,\ldots,p+1,
\end{eqnarray*}
with
\begin{eqnarray*}
\mu_{ir}^{*}=\int_0^{\tau}
\frac{\pi_r(s,\hat{\bgamma},\hat{\Lambda}_0) \hat{p}(s-)}
{\mathcal{Y}(s,\hat{\Lambda}_0)} \sum_{j=1}^{m_i} d\hat{M}_{ij}(s)
\end{eqnarray*}
and
\begin{eqnarray*}
\hat{M}_{ij}(t)=N_{ij}(t)-\int_0^{t}\exp({\hat{\bbeta}}^T
{\bf Z}_{ij})Y_{ij}(u)
\psi_i(\hat{\bgamma},\hat{\Lambda}_0,u-)d\hat{\Lambda}_0(u).
\end{eqnarray*}

\bsh
\underline{Step IV}
\esh

First order Taylor expansion of
$\bU(\hat{\bgamma},\hat{\Lambda}_0(\cdot,\hat{\bgamma}))$ about
$\bgamma^\circ=({\bbeta^\circ}^T,\theta^\circ)^T$ gives
 \begin{eqnarray*}
\bU(\hat{\bgamma},\hat{\Lambda}_0(\cdot,\hat{\bgamma}))
=\bU(\bgamma^\circ,\hat{\Lambda}_0(\cdot,\bgamma^\circ))+ \bD(\bgamma^\circ)
(\hat{\bgamma}-\bgamma^\circ)^T + o_p(1),
 \end{eqnarray*}
where
 \begin{eqnarray*}
D_{ls}(\bgamma)=\partial
U_l(\bgamma,\hat{\Lambda}_0(\cdot,\bgamma)) / \partial \gamma_s
 \end{eqnarray*}
for $l,s=1,\ldots,p+1$, with $\gamma_{p+1}=\theta$.

For $l,s=1,\ldots,p$ we have
 \begin{eqnarray}\label{dls}
\hspace*{-2em} D_{ls}(\bgamma)&=&-n^{-1}\sum_{i=1}^{n}
 \left\{
 \frac{\phi_{2i}(\bgamma,\hat{\Lambda}_0,\tau)}{\phi_{1i}(\bgamma,\hat{\Lambda}_0,\tau)}
 \sum_{j=1}^{m_i}Z_{ijl} \frac{\partial \hat{H}_{ij}(T_{ij})}{\partial
 \beta_s} \right. \nonumber\\
& & - \left. \left[
 \frac{\phi_{3i}(\bgamma,\hat{\Lambda}_0,\tau)}{\phi_{1i}(\bgamma,\hat{\Lambda}_0,\tau)}
 -\frac{\phi^2_{2i}(\bgamma,\hat{\Lambda}_0,\tau)}{\phi^2_{1i}(\bgamma,\hat{\Lambda}_0,\tau)}
 \right]
 \sum_{j=1}^{m_i} \hat{H}_{ij}(T_{ij}) Z_{ijl}\frac{\partial \hat{H}_{i.}(\tau)}{\partial
 \beta_s}
 \right\},
 \end{eqnarray}
 \begin{eqnarray*}
\frac{\partial \hat{H}_{ij}(\tau_k)}{\partial \beta_s} =
\frac{\partial \hat{\Lambda}_0(T_{ij} \wedge \tau_k)}{\partial
\beta_s} \exp(\bbeta^T {\bf Z}_{ij}) +\hat{\Lambda}_0(T_{ij}
\wedge \tau_k)\exp(\bbeta^T {\bf Z}_{ij}) Z_{ijs}
 \end{eqnarray*}
and
\begin{eqnarray*}
\frac{\partial \Delta \hat{\Lambda}_0(\tau_k)}{\partial \beta_s}
&=& -d_k\left\{\sum_{i=1}^{n}
 \frac{\phi_{2i}(\bgamma,\hat{\Lambda}_0,\tau_{k-1})}
 {\phi_{1i}(\bgamma,\hat{\Lambda}_0,\tau_{k-1})} R_{i.}(\tau_k)
 \right\}^{-2} \nonumber \\
&& \sum_{i=1}^n \left[ \left\{
 \frac{\phi_{2i}^2(\bgamma,\hat{\Lambda}_0,\tau_{k-1})}
 {\phi_{1i}^2(\bgamma,\hat{\Lambda}_0,\tau_{k-1})}
 -
 \frac{\phi_{3i}(\bgamma,\hat{\Lambda}_0,\tau_{k-1})}
 {\phi_{1i}(\bgamma,\hat{\Lambda}_0,\tau_{k-1})}
 \right\} \frac{\partial \hat{H}_{i.}(\tau_{k-1})}{\partial \beta_s} R_{i.}(\tau_k) \right. \nonumber \\
&& + \left.
 \frac{\phi_{2i}(\bgamma,\hat{\Lambda}_0,\tau_{k-1})}
 {\phi_{1i}(\bgamma,\hat{\Lambda}_0,\tau_{k-1})}
 \sum_{j=1}^{m_i} R_{ij}(\tau_{k}) Z_{ijs}
 \right].
\end{eqnarray*}
For $l=1,\ldots,p$ we have
\begin{eqnarray}\label{dlp}
 D_{l(p+1)}(\bgamma) &=&
 -n^{-1}\sum_{i=1}^{n}\left\{
 \frac{\phi_{2i}(\bgamma,\hat{\Lambda}_0,\tau)}{\phi_{1i}(\bgamma,\hat{\Lambda}_0,\tau)}
 \sum_{j=1}^{m_i} Z_{ijl} \frac{\partial \hat{H}_{ij}(T_{ij})}{\partial \theta} \right. \nonumber\\
 &&  \left.
 + \left[ \frac{\phi_{2i}^{(\theta)}(\bgamma,\hat{\Lambda}_0,\tau)}{\phi_{1i}(\bgamma,\hat{\Lambda}_0,\tau)}
 -\frac{\phi_{2i}(\bgamma,\hat{\Lambda}_0,\tau)\phi_{1i}^{(\theta)}(\bgamma,\hat{\Lambda}_0,\tau)}
  {\phi_{1i}^2(\bgamma,\hat{\Lambda}_0,\tau)} \right. \right. \nonumber\\
 && \left. \left.
  +\left\{\frac{\phi_{2i}^2(\bgamma,\hat{\Lambda}_0,\tau)}{\phi_{1i}^2(\bgamma,\hat{\Lambda}_0,\tau)}
  - \frac{\phi_{3i}(\bgamma,\hat{\Lambda}_0,\tau)}{\phi_{1i}(\bgamma,\hat{\Lambda}_0,\tau)} \right\}
  \frac{\partial \hat{H}_{i.}(\tau)}{\partial \theta} \right]
  \sum_{j=1}^{m_i}\hat{H}_{ij}(T_{ij})Z_{ijl}
 \right\}
\end{eqnarray}
and
 \begin{eqnarray}\label{dpl}
 D_{(p+1)l}(\bgamma) &=& n^{-1}\sum_{i=1}^{n}
 \left\{
 \frac{
 \phi_{1i}^{(\theta)}(\bgamma,\hat{\Lambda}_0,\tau) \phi_{2i}(\bgamma,\hat{\Lambda}_0,\tau)}
 {\phi_{1i}^2(\bgamma,\hat{\Lambda}_0,\tau)}
-
 \frac{
 \phi_{2i}^{(\theta)}(\bgamma,\hat{\Lambda}_0,\tau)}
 {\phi_{1i}(\bgamma,\hat{\Lambda}_0,\tau)}
 \right\} \frac{\partial \hat{H}_{i.}(\tau)}{\partial \beta_l}.
 \end{eqnarray}
Finally,
 \begin{eqnarray}\label{dpp}
 D_{(p+1)(p+1)}(\bgamma)&=&n^{-1}\sum_{i=1}^{n}
 \left\{
 \frac{
 \phi_{1i}^{(\theta,\theta)}(\bgamma,\hat{\Lambda}_0,\tau)
 }{\phi_{1i}(\bgamma,\hat{\Lambda}_0,\tau)}
  -
\left[ \frac{
 \phi_{1i}^{(\theta)}(\bgamma,\hat{\Lambda}_0,\tau)
 }{\phi_{1i}(\bgamma,\hat{\Lambda}_0)}\right]^2 \right. \nonumber \\
&& \left. +\left[ \frac{
\phi_{1i}^{(\theta)}(\bgamma,\hat{\Lambda}_0,\tau)
\phi_{2i}(\bgamma,\hat{\Lambda}_0,\tau)}{
\phi_{1i}^{2}(\bgamma,\hat{\Lambda}_0,\tau)}
 -\frac{\phi_{2i}^{(\theta)}(\bgamma,\hat{\Lambda}_0,\tau)}{\phi_{1i}(\bgamma,\hat{\Lambda}_0,\tau)} \right]
 \frac{\partial \hat{H}_{i.}(\tau)}{\partial \theta}
 \right\}
 \end{eqnarray}
where
 $$
\phi_{1i}^{(\theta,\theta)}(\bgamma,\hat{\Lambda}_0,\tau) =\int
w^{N_{i.}(\tau)}\exp\{-w \hat{H}_{i.}(\tau)\} \frac{d^2 f(w)}
 {d \theta^2} dw,
 $$
 \begin{eqnarray*}
\frac{\partial \hat{H}_{ij}(\tau_k)}{\partial \theta} =
\frac{\partial \hat{\Lambda}_0(T_{ij} \wedge \tau_{k})}{\partial
\theta} \exp(\bbeta^T {\bf Z}_{ij}),
 \end{eqnarray*}
and
\begin{eqnarray*}
\frac{\partial \Delta \hat{\Lambda}_0(\tau_k)}{\partial \theta}
&=& -d_k\left\{\sum_{i=1}^{n}
 \frac{\phi_{2i}(\bgamma,\hat{\Lambda}_0,\tau_{k-1})}
 {\phi_{1i}(\bgamma,\hat{\Lambda}_0,\tau_{k-1})} R_{i.}(\tau_k)
 \right\}^{-2} \nonumber \\
&& \sum_{i=1}^n
 R_{i.}(\tau_k)
 \left[
 \frac{\phi_{2i}^{(\theta)}(\bgamma,\hat{\Lambda}_0,\tau_{k-1})}
 {\phi_{1i}(\bgamma,\hat{\Lambda}_0,\tau_{k-1})}
 -
 \frac{\phi_{2i}(\bgamma,\hat{\Lambda}_0,\tau_{k-1})\phi_{1i}^{(\theta)}(\bgamma,\hat{\Lambda}_0,\tau_{k-1})}
 {\phi_{1i}^2(\bgamma,\hat{\Lambda}_0,\tau_{k-1})}
  \right.  \nonumber \\
 && \left.
 + \frac{\partial \hat{H}_{i.}(\tau_{k-1})}{\partial \theta}
 \left\{
  \frac{\phi_{2i}^2(\bgamma,\hat{\Lambda}_0,\tau_{k-1})}{\phi_{1i}^2(\bgamma,\hat{\Lambda}_0,\tau_{k-1})}
  - \frac{\phi_{3i}(\bgamma,\hat{\Lambda}_0,\tau_{k-1})}{\phi_{1i}(\bgamma,\hat{\Lambda}_0,\tau_{k-1})}
  \right\} \right].
\end{eqnarray*}

\bsh
\underline{Step V}
\esh

Combining the results above we get that
$n^{1/2}(\hat{\bgamma}-\bgamma^\circ)$ is asymptotically zero-mean
normally distributed with a covariance matrix that can be
consistently estimated by
 \begin{eqnarray*}
  \hat{\bD}^{-1}(\hat{\bgamma})\{\hat{\bV}(\hat{\bgamma})+\hat{\bG}(\hat{\bgamma})+\hat{\bC}(\hat{\bgamma})\}
  \hat{\bD}^{{-1}}(\hat{\bgamma})^T.
 \end{eqnarray*}

\subsection{Proof of (\ref{acn})}

The goal is to prove that
\begin{equation}
\sup_{s,\bgams} |A(\bgamma,\tilde{\Lambda}^{(n)},s) - a(\bgamma,\tilde{\Lambda}^{(n)},s)| \rightarrow 0
\hspace*{1em} \mbox{a.s.}
\label{acnn}
\end{equation}
This involves several steps.

First, it is easy to see that there exists a constant $\kappa$ (independent of
$\bgamma$ and $s$) such that
\begin{eqnarray}
\sup_{s,\bgams} | A(\bgamma,\Lambda_1,s) - A(\bgamma,\Lambda_2,s) |
& \leq & \kappa \| \Lambda_1 - \Lambda_2 \|, \label{Alip} \\
\sup_{s,\bgams} | a(\bgamma,\Lambda_1,s) - a(\bgamma,\Lambda_2,s) |
& \leq & \kappa \| \Lambda_1 - \Lambda_2 \|. \label{alip}
\end{eqnarray}
Next, for any fixed continuous $\Lambda$, the functional strong law of large numbers of
Andersen \& Gill (1982, Appendix III) implies that, with probability one,
\begin{equation}
\sup_{s,\bgams} |A(\bgamma,\Lambda,s) - a(\bgamma,\Lambda,s)| \rightarrow 0.
\label{acc}
\end{equation}

Now, given $\eps>0$, define the sets $\{ t_j^\peps \}$, $\{ \bgamma_k^\peps \}$,
and $\{ \Lambda_l^\peps \}$
to be finite partition grids of $[0,\tau]$, $\cG$, and
$[0,\Lambda_{max}]$, respectively, with distance of no more than $\eps$ between grid
points. Define $\cL_\eps^*$ to be the set of functions of $t$ and $\bgamma$ defined by
linear interpolation through vertices of the form $(t_j^\peps,\bgamma_k^\peps,\Lambda_l^\peps)$.

Obviously $\cL_{\eps}^*$ is a finite set. Hence, in view of (\ref{acc}), there exists
a probability-one set of realizations $\Omega_\eps$ for which
\begin{equation}
\sup_{s \in [0,\tau], \bgams \in \cG, \Lambda \in \cL_\eps^*}
|A(\bgamma,\Lambda,s) - a(\bgamma,\Lambda,s)| \rightarrow 0.
\label{acle}
\end{equation}
Define
$$
\Omega^{**} =
\bigcap_{\ell=1}^\infty \Omega_{1/\ell}
$$
and $\Omega_0 = \Omega^* \cap \Omega^{**}$, with $\Omega^*$ as defined earlier.
Clearly $\Pr(\Omega_0)=1$. From now on, we restrict attention to $\Omega_0$.

Now let $\eps>0$ be given. Choose $\ell > \eps^{-1}$. In view of (\ref{equi}) and
(\ref{acle}), we can find for any $\omega \in \Omega_0$ a suitable positive integer
$\bar{n}(\eps,\omega)$ such that, whenever $n \geq \bar{n}(\eps,\omega)$,
\begin{equation}
|\tilde{\Lambda}_0^{(n)}(t,\bgamma) - \tilde{\Lambda}_0^{(n)}(u,\bgamma)|
\leq B^* (t-u) + \frac{\epsilon}{2} \quad \forall t,u,
\label{thingone}
\end{equation}
\begin{equation}
\sup_{s \in [0,\tau], \bgams \in \cG, \Lambda \in \cL_{1/\ell}^*}
|A(\bgamma,\Lambda,s) - a(\bgamma,\Lambda,s)| \leq \eps.
\label{thingtwo}
\end{equation}

Next, let $\bar{\Lambda}_0^{(n)}$ denote the function defined by linear interpolation
through $(t_j^\peps,\bgamma_k^\peps,\bar{\Lambda}_{jk}^\peps)$, where $\bar{\Lambda}_{jk}^\peps$
is the element of $\{ \Lambda_l^\peps \}$ that is closest to $\tilde{\Lambda}_0^{(n)}
(t_j^\peps,\bgamma_k^\peps)$. It is clear that
$$
|\bar{\Lambda}_0^{(n)}(t_j^\peps,\bgamma_k^\peps)
- \tilde{\Lambda}_0^{(n)}(t_j^\peps,\bgamma_k^\peps)| \leq \eps \quad \forall j,k.
$$
Using (\ref{thingone}) and the Lipschitz continuity of
$\tilde{\Lambda}_0^{(n)}(t,\bgamma)$ with respect to $\bgamma$
(which follows from the corresponding property of $\hat{\Lambda}_0(t,\bgamma))$,
we thus obtain
$$
\sup_{t,\bgams} |\bar{\Lambda}_0^{(n)}(t,\bgamma) - \tilde{\Lambda}_0^{(n)}(t,\bgamma)|
\leq B^{**} \eps
$$
for a suitable fixed constant $B^{**}$ (depending on $B^*$ and $C^*$).
Combining this with (\ref{thingtwo}) and (\ref{alip}), we obtain
$$
\sup_{s,\bgams} |A(\bgamma,\tilde{\Lambda}^{(n)},s) - a(\bgamma,\tilde{\Lambda}^{(n)},s)|
\leq (2\kappa B^{**} + 1) \eps \quad \mbox{for all } n \geq \bar{n}(\eps,\omega).
$$
Since $\eps$ was arbitrary, the desired conclusion (\ref{acnn}) follows, and the proof is
thus complete.

\subsection{Definition and behavior of $h_i(r,s)$}

The quantity $h_i(r,s)$ appearing in (\ref{eq:yapp}) is given by
\begin{eqnarray*}
h_i(r,s)&=&\frac{2R_{i.}(s)\eta_{1i}(r,s)}{\{{\mathcal
Y}(s,\Lambda_0^\circ+r\Delta)\}^3}\frac{1}{n}\sum_{l=1}^{n}R_{l.}(s)\eta_{1l}(r,s)\sum_{j=1}^{m_i}
\exp(\bbeta^T {\bf Z}_{lj})\Delta(T_{lj}\wedge s) \nonumber \\
 && - \frac{R_{i.}(s)\eta_{2i}(r,s)}{\{{\mathcal
Y}(s,\Lambda_0^\circ+r\Delta)\}^2}\sum_{j=1}^{m_i} \exp(\bbeta^T {\bf
Z}_{ij})\Delta(T_{ij}\wedge s)
\end{eqnarray*}
where $\Delta(T_{ij}\wedge s)=\hat{\Lambda}_0(T_{ij}\wedge
s)-\Lambda_0^o(T_{ij}\wedge s)$ and
 $$
\eta_{2i}(r,s)=2\left\{\frac{\phi_{2i}(\bgamma^\circ,{\Lambda}_0^\circ+r\Delta,s)}
{\phi_{1i}(\bgamma^\circ,{\Lambda}_0^\circ+r\Delta,s)}\right\}^3
+\frac{\phi_{4i}(\bgamma^\circ,{\Lambda}_0^\circ+r\Delta,s)}
{\phi_{1i}(\bgamma^\circ,{\Lambda}_0^\circ+r\Delta,s)}
-3\frac{\phi_{2i}(\bgamma^\circ,{\Lambda}_0^\circ+r\Delta,s)\phi_{3i}(\bgamma^\circ,{\Lambda}_0^\circ+r\Delta,s)}
{\left\{\phi_{1i}(\bgamma^\circ,{\Lambda}_0^\circ+r\Delta,s)\right\}^2}.
 $$
For all $i=1,\ldots, n$ and $s \in [0,\tau]$,
we have $0 \leq R_{i.}(s) \leq m \nu$,
where $\nu$ is as in (\ref{eb}).
Moreover, for $k=1,\ldots,4$, we have
 $$
 \E[W_i^{r_{min}+(k-1)}\exp\{-W_i m e^{\sbeta^T Z} \Lambda_0^\circ(\tau)\}]
 \leq \phi_{ki}(\bgamma^\circ,{\Lambda}_0^\circ,s) \leq
 \E[W_i^{r_{max}+(k-1)}]
  $$
where $r_{\max}=\arg \max_{1\leq r \leq m} \E (W_i^r)$,
$r_{\min}=\arg \min_{1\leq r \leq m} \E (W_i^r)$. Hence,
$\eta_{1i}$ and $\eta_{2i}$ are bounded. In addition, the the
proof of Lemma 2 show that ${\mathcal Y}(s,\Lambda^\circ
+r\Delta)$ is uniformly bounded away from zero for $n$
sufficiently large. Finally, in the consistency proof we obtained
$\| \Delta \| = o(1)$. Therefore $h_{i}(r,s)$ is $o(1)$ uniformly
in $r$ and $s$.

\section{References}
\begin{description}
 \item
  {\sc Aalen, O. O.} (1976). Nonparametric inference in
connection with multiple decrement models. {\em Scand. J.
Statist.} {\bf 3}, 15-27.
  \item
  {\sc Abramowitz, M. and Stegun,
I. A. (Eds.)} (1972). {\em Handbook of Mathematical Functions with
Formulas, Graphs, and Mathematical Tables, 9th printing} New York:
Dover.
 \item
  {\sc Andersen, P. K., Borgan, O, Gill, R. D. and
Keiding, N.} (1993). Statistical models based on counting
processes. Berlin: Springer-Verlag.
 \item
  {\sc Andersen, P. K. and
Gill, R. D.} (1982). Cox's regression model for counting
processes: A large sample study. {\em  Ann. Statist.} {\bf 10},
1100-1120.
 \item {\sc Andersen, P. K., Klein, J. P., Knudsen, K.
M. and Palacios, R. T.} (1997). Estimation of variance in Cox's
regression model with shared gamma frailty. {\em Biometrics} {\bf
53}, 1475-1484.
 \item {\sc Breslow, N.} (1974). Covariance
analysis of censored survival data. {\em Biometrics,} {\bf 30},
89-99.
 \item {\sc Fine, J. P., Glidden D. V. and Lee, K.} (2003).
A simple estimator for a shared frailty regression model. {\em J.
R. Statist. Soc.} {\bf B 65}, 317-329.
 \item {\sc Cox, D. R.}
(1972). Regression models and life tables (with discussion). {\em
J. R. Statist. Soc.} {\bf B 34}, 187-220.
 \item {\sc Foutz, R. V.}
(1977). On the unique consistent solution to the likelihood
equation. {\em J. Amer. Statist. Ass.} {\bf 72}, 147-148.
 \item
{\sc Gill, R. D.} (1985). Discussion of the paper by D. Clayton
and J. Cuzick. {\em J. R. Statist. Soc.} {\bf A 148}, 108-109.
 \item {\sc Gill, R. D.} (1989). Non- and semi-parametric maximum
likelihood estimators and the Von Mises method (Part 1). {\em
Scand. J. Statist.} {\bf 16}, 97-128.
 \item {\sc Gill, R. D.}
(1992). Marginal partial likelihood. {\em Scand. J. Statist.} {\bf
79}, 133-137.
 \item {\sc Gorfine, M., Zucker, D. M., and Hsu, L.}
(2006). Prospective survival analysis with a general
semiparametric shared frailty model - a pseudo full likelihood
approach. {\em Biometrika}, to appear.
 \item {\sc Hartman, P.}
(1973). {\em Ordinary Differential Equations,} 2nd ed. (reprinted,
1982), Boston: Birkhauser.
 \item {\sc Henderson, R. and Oman, P.}
(1999). Effect of frailty on marginal regression estimates in
survival analysis. {\em J. R. Statist. Soc.} {\bf B 61}, 367-379.
\item
{\sc Hougaard, P.} (1986). Survival models for heterogeneous
populations derived from stable distributions. {\em Biometrika}
{\bf 73}, 387-396.
\item
{\sc Hougaard, P.} (2000). {\em Analysis of Multivariate Survival
data}. New York: Springer.
\item
{\sc Klein, J. P.} (1992). Semiparametric estimation of random
effects using the Cox model based on the EM Algorithm. {\em
Biometrics} {\bf 48}, 795-806.
\item
{\sc Louis, T. A.} (1982). Finding the observed information matrix
when using the EM algorithm. {\em J. R. Statis. Soc.} {\bf B 44},
226-233.
\item
{\sc McGlichrist, C. A.} (1993). REML estimation for survival
models with frailty. {\em Biometrics} {\bf 49}, 221-225.
\item
{\sc Murphy, S. A.} (1994). Consistency in a proportional hazards
model incorporating a random effect. {\em Ann. Statist.} {\bf 22},
712-731.
\item
{\sc Murphy, S. A.} (1995). Asymptotic theory for the frailty
model. {\em Ann. Statist.} {\bf 23}, 182-198.
\item
{\sc Nielsen, G. G., Gill, R. D., Andersen, P. K. and Sorensen, T.
I.} (1992). A counting process approach to maximum likelihood
estimation of frailty models. {\em Scand. J. Statist.} {\bf 19},
25-43.
\item
{\sc Parner, E.} (1998). Asymptotic theory for the correlated
gamma-frailty model. {\em Ann. Statist.} {\bf 26}, 183-214.
\item
{\sc Ripatti, S. and Palmgren J.} (2000). Estimation of
multivariate frailty models using penalized partial likelihood.
{\em Biometrics} {\bf 56}, 1016-1022.
\item
{\sc Vaida, F. and Xu, R. H.} (2000). Proportional hazards model
with random effects. {\em Stat. in Med.} {\bf 19}, 3309-3324.
\item
{\sc Yang, S. and Prentice, R. L.} (1999). Semiparametric
inference in the proportional odds regression model. {\em J. Amer.
Statist. Ass.} {\bf 94}, 125-136.
\item
{\sc Zucker, D. M.} (2005). A pseudo partial likelihood method for
semi-parametric survival regression with covariate errors. {\em
J. Amer. Statist. Ass.} {\bf 100}, 1264-1277.
\end{description}

\vspace*{0.5in}

\begin{flushleft}
David M. Zucker \\
Department of Statistics \\
Hebrew University \\
Mt. Scopus, 91905 Jerusalem, ISRAEL \\
E-mail: mszucker@mscc.huji.ac.il
\end{flushleft}

\end{document}